\documentclass{amsart}
\usepackage[english]{babel}
\usepackage{amsfonts, amssymb}
\usepackage[mathcal]{eucal}

\theoremstyle{plain}
\newtheorem{thm}{Theorem}

\theoremstyle{definition}

\newcommand\N{{\mathbb N}}

\newcommand\M{{\mathfrak M}}
\newcommand\cs{{\mathcal S}}
\newcommand\s{{\mathfrak S}}

\newcommand\bx{{\bf X}}

\newcommand\eset{\varnothing}
\title{A non constructive proof of the existence of a maximal irreducibility
measure.}
\author{Sanatan Rai}
\email{sanat@stanford.edu}
\address{Department of Management Science and Engg.\\%
Stanford University\\%
Stanford, California 94305.}

\begin{document}

\begin{abstract}
We deduce the existence of a maximal irreducibility measure for a
Markov chain using Zorn's lemma.
\end{abstract}

\maketitle

Let $\bx :=\{X_n\}_{n\in \N_0}$ be a Markov Chain with values in 
space $(\cs,\s)$. A measure $\phi$ on $(\cs,\s)$ is called an
\emph{irreducibility measure} for $\bx$, if for all $A\in\s$,
\begin{equation}
\phi(A)>0\implies (\forall x\in\cs)\ L(x,A)>0.\label{eq:irred}
\end{equation}
Here $L(x,A):=P_x\{\tau_A<\infty\}$, for $\tau_A:=\inf\{n\ge 0:X_n\in A\}$. 
In general, there will be many
irrducibility measures for the given chain. These measures tell us
about the structure of the space $(\cs, \s)$ as seen by $\bx$. So which
measure is the best? 

We have not said anything about the reverse implication in (\ref{eq:irred}). It
is natural to question when the reverse implication is true. The reverse implication
is true for the \emph{maximal irreducibility measure}, and is in fact a characterising
feature. So when does such a measure exist? As long as there is \emph{at least one}
irreducibility measure for $\bx$, there is a maximal irreducibility measure. 

Note that if $\mu \ll \psi$ is another measure, then $\mu$ is also an irreducibility
measure. This is so  because, if $\mu(A)>0$ then $\phi(A)>0$, but then $L(x,A)>0$ for
all $x$.  We can now state the (\cite{meyntweed})

\begin{thm} If there exists an irreducibility measure $\phi$ on $\s$ for the
chain $\bx$, then there exists a unique maximal irreducibility measure $\psi$ on 
$\s$ such that
\begin{enumerate}
\item $\psi$ is an irreducibility measure for $\bx$, 
\item $\mu$ is an irreducibility measure for $\bx$ if and only if $\mu \ll \psi$, and,
\item if $A\in\s$ is such that $\psi(A)=0$, then $\psi\{x:L(x,A)>0\}=0$.
\end{enumerate}
\end{thm}
This theorem also defines a maximal irreducibility measure. 

\begin{proof} The uniqueness follows from (2) above. For if $\psi$ and $\psi'$ are maximal
irreducibility measures, $\psi\ll\psi'$ and $\psi'\ll\psi$, since they are also 
irreducibility measures. Hence $\psi\sim\psi'$. 

Let $\M$ be the set of all irreducibility measures on $\s$ for $\bx$. Impose the
order $\prec$ ($\succ$) by $\phi\prec\phi'\equiv\phi'\succ\phi\Leftrightarrow\phi\ll\phi'$.

\noindent{\sc Claim:} Every chain has an upper bound.
\begin{proof}[Proof of claim] Let $\phi_1\prec\phi_2\prec\ldots$ be a chain. 
Define the upper bound $\phi$ as follows.
\begin{equation}
\forall A\in\s,\ \phi(A):=\lim_n\phi_n(A)\label{eq:phimax}
\end{equation}
If we can show that $\phi$ is a measure on $\s$, then $\phi$ will be the required upper
bound. For the $\phi_m$ form a chain, $\phi_m(A)\le\phi_{m+1}(A)$ and so
the sequence of numbers $\phi_m(A)$ is increasing, and so, $\phi(A)=\lim_k\phi_k(A)
\ge\phi_m(A)$ for all $m\in\N$. Hence, if $\phi(A)=$, then $\phi_m(A)=$ for all $m$. This
implies that $\phi\gg \phi_m,\ \forall m\in\N$.  

To show that $\phi$ is a measure, we must prove that
\begin{enumerate}
\item $\phi(\eset)=0$,
\item $A,B\in\s:A\subset B\implies\phi(A)\le\phi(B)$, 
\item and, $\{A_n\}_{n\in\N}\subseteq\s$ then $\phi(\cup_n A_n)\le\sum_n\phi(A_n)$
\end{enumerate}

To this end, note that $\phi(\eset)=\lim_n\phi_n(\eset)=\lim_n 0 = 0$, and for 
$A\subset B$, $\phi(A)=\lim_n(A) \le \lim_n\phi_n(B)=\phi(B)$. Similarly,
$\phi(\cup_k A_k)=\lim_n\phi_n(\cup_k A_k)\le\lim_n\sum_k\phi_n(A_k)=\sum_k\lim_n\phi_n(A_k)
=\sum_k\phi(A_k)$, where the interchange of the sum and the limit is allowed since all
the terms are positive. 

This shows that $\phi$ is indeed a measure.
\end{proof}

Thus we have shown that every chain has an upper bound. Since there exists \emph{some} 
irreducibility measure, the set of elements ordered by $\gg$ is well not empty. So by Zorn's 
lemma, there exists a maximal element. This element is the maximal irreducibility measure. 

\end{proof}


\begin{thebibliography}{1}
\bibitem{meyntweed} Meyn, Sean P. and R.~ L.~ Tweedie: {\it Markov
Chains and Stochastic Stability}, Springer, 1993. 
\end{thebibliography}
\end{document}